\newcommand{\C}{\mathbb C}

\newcommand{\R}{\mathbb R}

\newcommand{\Z}{\mathbb Z}

\newcommand{\si}{\sigma}
\newcommand{\la}{\lambda}

\newcommand{\Ind}{\mathrm{Ind\,}}
\newcommand{\J}{J\,}

\def\01{[0,1]}
\def\!{^{-1}}

\newcommand{\cf}{\mathcal{F}}

\newcommand{\cg}{\mathcal{G}}
\newcommand{\cL}{\mathcal{L}}
\newcommand{\cLL}{\mathcal{L_\lambda}}

\newcommand{\cB}{\mathcal{B}}

\newcommand{\cS}{\mathcal{S}}

\newcommand{\hk}{H^{k+s}(\Omega;\C^m )}

\newcommand{\im}{\operatorname{Im}}

\newcommand{\hfor}{\hbox{ \,for\, }}
\newcommand{\hand}{\hbox{ \,and\, }}

\newcommand{\LL}{\Lambda}
\newcommand{\tr}{\operatorname{tr}}

\newcommand{\pc}{\,Ell(\R^n)\,}
\newcommand{\ch}{\operatorname{ch}}

\newcommand{\bt}{\begin{theorem}}
\newcommand{\et}{\end{theorem}}
\newcommand{\bc}{\begin{corollary}}
\newcommand{\ec}{\end{corollary}}
\newcommand{\bp}{\begin{proposition}}
\newcommand{\ep}{\end{proposition}}
\newcommand{\bl}{\begin{lemma}}
\newcommand{\el}{\end{lemma}}
\newcommand{\br}{\begin{remark}}
\newcommand{\er}{\end{remark}}
\newcommand{\bd}{\begin{definition}}
\newcommand{\ed}{\end{definition}}
\newcommand{\be}{\begin{equation}}
\newcommand{\ee}{\end{equation}}
\newcommand{\ra}{\rightarrow}

 \documentclass[thmsa, a4paper]{amsart}
 
\usepackage{amsmath}
\usepackage{amsfonts}
\usepackage{amssymb}
\input{diagxy}

\begin{document}
\newtheorem{theorem}{Theorem}[section]
\newtheorem{corollary}[theorem]{Corollary}
\newtheorem{lemma}[theorem]{Lemma}
\newtheorem{proposition}[theorem]{Proposition}
\newtheorem{remark}{Remark}[section]
\newtheorem{definition}{Definition}[section]
\numberwithin{equation}{section}
\numberwithin{theorem}{subsection}
\numberwithin{definition}{subsection}
\numberwithin{remark}{subsection}

\title[Index Bundle and Bifurcation II]{Bifurcation of Fredholm Maps II; The  Dimension of the Set of Bifurcation Points}

\author{J. Pejsachowicz }

\address{Dipartimento di Matematica Politecnico di Torino\\
    Corso Duca degli Abruzzi  24\\10129 Torino, Italy.}
\email{jacobo.pejsachowicz@polito.it}
\thanks{This work was supported by  MIUR-PRIN2007-Metodi variazionali e topologici nei fenomeni nonlineari}


\date{\today}
\begin{abstract}
We obtain an estimate  for the covering dimension of the set of bifurcation points  for  solutions of  nonlinear elliptic boundary value problems  from the principal symbol of the linearization along the trivial branch of solutions.
\end{abstract}

\maketitle

\section{Introduction}

 In \cite{[Pe]} we defined an index of bifurcation points  $\beta(f)$ of a  parametrized family $f$  of $C^1$-Fredholm maps.  Nonvanishing of  $\beta(f)$ entails the existence of at least one point of bifurcation from a trivial branch of zeroes of the family  $f.$ Linearization of  $f$ along  the trivial branch produces a parametrized family of linear Fredholm operators 
 $L.$ The  index  $\beta(f)$  depends only on the  stably  fiberwise homotopy equivalence class of the index bundle $\Ind L$ of $L.$  The nonvanishing of $\beta(f)$ can be checked  through  the Stiefel-Whitney  and Wu  characteristic classes of  $\Ind L$ since they  are invariant under stably fiberwise homotopy equivalence.  
 
If the parameter space  is a manifold, the nonvanishing  of the  Stiefel-Whitney and Wu classes of $\Ind L$ not only implies that the set $B(f)$ of all  bifurcation points of $f$ is nonempty, but also provides some further information about the covering dimension of this set and its position in the parameter space.  
 
In this paper, using  the above observation  (Theorem \ref{th:1.2} in section 3)  together with our results from \cite{[Pe]}, we obtain an estimate on the covering dimension of the set of bifurcation points of solutions of nonlinear elliptic boundary value problems parametrized by a smooth manifold.
 
  Roughly speaking, the approach is as follows:  assuming that  the coefficients of the leading terms of $L$  are independent from the parameter near the boundary, an extension to families of the Agranovich reduction identifies the complexification of $\Ind L$  with the index bundle of a family of pseudo-differential operators $S$ on $\R^n.$  Applying to $S$ a cohomological form of the Atiyah-Singer family index theorem, due to Fedosov,  we determine the  Chern character of the complexification $c(\Ind L) $ as an integral along the fiber of a differential form associated to the symbol of the family. 
  
   In principle, the above  leads to the computation of Wu classes of the index bundle of  $L,$  since they are polynomials in Pontriagin classes of $\Ind L$ with $Z_{p}$ coefficients.  However, the general  expression is messy and can hardly be used in practice.  It becomes much simpler  by evaluating  Wu classes of $\Ind L$ on spherical homology classes. Restricting the  family $f$ to spheres embedded in the parameter space   and using  our approach in \cite{[Pe]}, we obtain explicit conditions for non vanishing of Wu classes  and hence estimates for the dimension of the set of bifurcation points.  A similar use of the Stiefel-Whitney classes gives some complementing results.

  The topological  dimension of the set of solutions of nonlinear equations  and  the set of bifurcation points has been discussed in various places, mainly in the case of  compact vector fields  and  semi-linear  Fredholm maps \cite{[Fi-Pe],[Ba],[Pe-1],[Iz-1],[Al-An],[Fi-Ma-Pe]}. However, it should be remarked, that our estimates are obtained directly  from the leading coefficients of  linearized equations without  the need to solve them.  This is the main reason for using elliptic invariants in a topological approach to bifurcation which complements  the classical Lyapunov-Schmidt method.

 The paper is organized as follows:   in section $2$  we state our main result, theorem  \ref{th:0}.  In section $3$  we relate  Wu classes of  $\Ind L$  to the  dimension of the set of bifurcation points of $f$.  In section $4,$  after discussing  the Agranovich reduction, we carry out the computations of the relevant  characteristic  classes,  completing in this way the proof of  theorem \ref{th:0}.  In section $5$, assuming  that the linearization along the trivial branch is a (real) lower-order  perturbation of a family of elliptic boundary value problems with complex coefficients, we  obtain sufficient conditions for bifurcation in  dimensions not covered by theorem \ref{th:0}, using  Stiefel-Whitney  classes. 
 
  I would like to thank Yuli Rudiak, Michael Crabb and Friedrich Hegenbarth  for their generous help.
  
\vskip 5pt

\section{The main theorem}

 We  consider boundary value problems of the form 
 \begin{equation}\label{bvp1}
\left\{\begin{array}{l} \cf\,(\la, x,u,\ldots ,D^{k}u)=0 \hfor  x \in  \Omega , \\ \cg^i(\la, x,u,\ldots, D^{k_i} u)=0 \hfor  x \in  \partial \Omega , \, 1 \le  i \le r,\end{array}\right. \end{equation} 
where $\Omega$ is  an open bounded subset of $\R^n$ with smooth boundary $\partial \Omega,$  $u\colon\bar\Omega\rightarrow \R^m$ is a vector function,  $\la$ is a parameter belonging to a smooth compact connected  $d$-dimensional manifold $\LL$ and, denoting with $k^*$ the number of multindices $\alpha$'s  such that $ |\alpha| \leq k,$  $$ \cf \colon \LL \times \bar \Omega \times  \R^{mk^*}  \rightarrow  \R^m \hand  \cg^{i} \colon \LL\times \bar\Omega  \times  \R^{mk_i^*} \rightarrow  \R  $$  are smooth with  $\cf (\la,x, 0) = 0,\  \cg^i( \la,x,0)=0, \, 1 \le  i \le r.$

Here and below we will freely use the notation  from \cite{[Pe]}.

\vskip 3pt  

We will denote with $(\cf,\cg)$  the family of nonlinear differential operators  \be\label{cinf} (\cf,\cg)\colon \R^{q}\times C^{\infty}(\bar\Omega; \R^{m}) \ra C^{\infty}(\bar \Omega; \R^{m})\times  C^{\infty}(\partial \Omega; \R^{r}) \ee defined by
\[(\cf,\cg)(\la, u)=( \cf(\la, x,u,..,D^{k}u), \tau \cg^1(\la, x,u,..,D^{k_{1}}u),.., \tau \cg^r(\la, x,u,..,D^{k_{r}}u),\] where $\tau$ is the restriction to the boundary.

 We assume:
\begin{itemize}

\item[$H_1)$]  
 For all $\la \in \Lambda,$  the linearization $(\mathcal L_\la(x,D), \mathcal B_\la(x,D))$  of $(\cf_\la, \cg_\la)$  at $u\equiv 0,$ is an elliptic boundary value problem.  Namely,   ${\mathcal L}_\la(x,D)$ is  elliptic, properly elliptic at the boundary, and the  boundary operator \[\mathcal B_\la(x,D)=(\mathcal B^1_\la(x,D),...,  \mathcal B^r_\la(x,D))^t\] verifies the Shapiro-Lopatinskij  condition with respect to $ {\mathcal L}_\la(x,D)$.

\item [$H_2)$]   There exists a point    $\nu \in \Lambda $ such that, for every  $f \in C^\infty(\bar \Omega; \C^m)$ and \\$g \in C^\infty (\partial\Omega; \C^r),$  the problem: 
\[\left\{ \begin{array}{l}{\mathcal L}_\nu (x,D) u(x) = f(x)\ \text{for}\  x \in \Omega \\ 
\mathcal {B}_\nu (x,D) u(x) = g(x)\ \text{for}\  x\in \partial\Omega,
\end{array}\right. \]
has  a unique smooth solution.
\item [$H_3)$]  \begin{itemize}
\item[ i)]  The coefficients  $b^{i}_\alpha(\la,x),  \  |\alpha|=k_i,  1\leq i\leq r ,$
 of the  leading terms of boundary operators $\mathcal B^1_\la(x,D),\dots, \mathcal B^r_\la(x,D)$ 
 are independent of  $\la.$  

 \item[ ii)]  There exist a compact set  $K\subset \Omega $ such that  the coefficients \\$a_{\alpha }(\la, x); \, |\alpha |=k $ of the  leading terms of $\cLL,$ are independent of $\la $ for $x\in \bar\Omega-K.$ \end{itemize}\end{itemize}

 Let $p(\la,x,\xi ) \equiv \sum^{}_{\vert \alpha \vert =k}a_{\alpha }(\la,x)\xi ^{\alpha }$ be the principal symbol of $\cLL.$  
  
  By ellipticity,  $ p(\la,x,\xi )\in GL(m;\C)$ if $\xi\neq 0.$ On the other hand, by $H_3,$
    \[ p (\la,x,\xi )= p (\nu,x,\xi )  \hfor x \in \bar\Omega-K.\]
      Therefore, putting   \[\sigma(\la,x,\xi)= Id  \ \text{ for any}\  (\la,x,\xi) \  \text{with } \  x\notin K, \]  the map  $\si(\la,x,\xi) = p(\la,x,\xi) p(\nu,x,\xi)^{-1}$ extends to a smooth map  
\begin{equation} \label{symbo}  \sigma \colon \LL\times (\R^{2n} - K\times\{0\}) \rightarrow GL(m;\C).\end{equation}

Assuming, without loss of generality,  that   $K\times\{0\}$ is contained in the unit ball $B^{2n} \subset \R^{2n},$ we associate to  $\sigma$   the restriction (pulback ) of the matrix one-form  $ \sigma^{-1} d \sigma$ to  $\Lambda \times \partial B^{2n} \simeq \Lambda\times S^{2n-1},$  which will  be denoted in the same way.
    
  Taking the trace of the  $(q+2n-1)$-th power of the matrix  $ \sigma^{-1} d\sigma$ we obtain an ordinary $(q+2n-1)$-form  $ \tr (\sigma^{-1}d\sigma)^{q+2n-1} $ on $\Lambda \times S^{2n-1}.$ 

Let $ \Sigma\simeq  S^{q} \subset \Lambda $ be an embedded sphere of even dimension  $q.$ We define the {\it  degree} $\deg(\sigma;\Sigma)$  of  $\si$ on $\Sigma$ by
   \begin{equation} \label{fed}
  \deg ( \sigma;\Sigma)=  \displaystyle {\frac{(\frac{1}{2}q+n-1)!} {(2\pi i)^{(\frac{1}{2}q+n)}( q+2n-1)!}} \int_{\Sigma \times S^{2n-1}} tr (\sigma^{-1}d \sigma)^{q+2n-1}.
  \end{equation}
  
We will see later  that, for any embedded sphere $\Sigma,$ $\deg( \sigma; \Sigma)$ is  an  integral number. 

\vskip 3pt

 Let us recall that a {\it bifurcation point from the trivial branch} for solutions of \eqref{bvp1} is a point $ \la_* \in \Lambda$ such that there exist a sequence  $(\la_n, u_n)\in\Lambda \times \C^\infty (\bar\Omega)$ of solutions  of  \eqref{bvp1}  with $u_n \neq 0,$  $\la_n \ra \la_*$ and $u_n \ra 0 $ uniformly with all of its derivatives. 

\begin{theorem}\label{th:0}
Let  the boundary value problem \eqref{bvp1} verify $H_1- H_3$ and let $p$ be an odd  prime such that  
$p\leq  d/2 +1.$ 

 If, $\LL$ is orientable and, for some embedded sphere  $\Sigma \subset \Lambda,$  of dimension\\ $q=2(p-1),$  
  $\deg(\sigma;\Sigma)$  is not divisible by $p,$   then
\begin{itemize}  
  \item[i)] \  the Lebesgue covering dimension of the set $B$ of 
all bifurcation points of \eqref{bvp1} is at least $d- q,$
\item[ii)]\ the set $B$ either disconnects $ \Lambda$ or is not contractible in  $ \Lambda$ to a point. 
\end{itemize}
\end{theorem}
\vskip 5pt
  
\section{Characteristic classes and bifurcation of Fredholm maps} \label{bfm}

We begin with a short recapitulation  of  \cite{[Pe]}. From now on, Fredholm means Fredholm of index $0.$ 

 Let $O$ be an open subset of a Banach space  $X$  and let  $ \Lambda$ be a  finite connected $CW$-complex. 
A family of $C^1$-Fredholm maps continuously parametrized by $\Lambda$ is a continuous map $f\colon \Lambda \times O \ra X$  such that the map  $f_{\la}\colon O\ra Y$ defined by $f_{\la}(x)= f(\la ,x) $ is  $C^1,$   for each ${\la} \in  \Lambda.$ Moreover,  $Df_{\la}(x)$ is a Fredholm operator of index $0$  which depends continuously on $(\la,x)$ with respect to the norm topology in the space of $\mathcal L (X,Y).$ 

We will  assume everywhere in this paper that  $O$ is  a neighborhood of the origin  $0\in X$ and  that $f(\la ,0)=0$  for all ${\la} $ in  $\Lambda.$ Solutions $(\la,0)$ of the  equation $f(\la ,x)=0$ form  the {\it trivial branch}, which we will identify  with  the parameter space  $\Lambda.$

  A point ${\la} _*$ in $ \Lambda$  is called   {\it bifurcation point from the trivial branch} for solutions of the equation $f(\la ,x)=0$ if  every neighborhood of $(\la _*,0)$ contains nontrivial solutions of this equation. 
  
 The   {\it  linearization of the family $f$  along the trivial branch } is the family of operators  $L\colon \Lambda\ra \Phi_0(X,Y)$ defined  by $L_\la =Df_{\la}(0),$ where the right hand side denotes  the Frechet derivative of $f_{\la}$ at $0.$ 
 
 Bifurcation  only occurs  at points ${\la} \in  \Lambda$  where  $L_{\la}$ is singular but, in general,  the set  $B(f)$ of all bifurcation points of  a family $f$ is only a proper  closed subset  of  the set $ \Sigma(L)$ of all singular points of  $L.$ 

 Given a compact  space $\Lambda,$  let  ${KO}(\Lambda)$ (resp. \!$K(\Lambda)$) be the Grothendieck group of all  real (resp. complex) virtual vector bundles over $\Lambda,$ and let $\tilde{KO}(\Lambda)$(resp. \!$\tilde K(\Lambda)$   be the corresponding  reduced group, i.e., the kernel of the rank homomorphism.   Recalling  that two  vector bundles are   {\it stably equivalent}  if  they become isomorphic after addition of trivial bundles to both sides, in this paper,  we will identify   $\tilde{KO}(\Lambda)$ with the group of all stable equivalence classes of vector bundles over $\Lambda.$

 With the above identification the {\it index bundle}  $\Ind L$ of a family $L$ of Fredholm operators  is defined as follows:  using  compactness of $ \Lambda,$ one can find a finite dimensional subspace $V$ of $Y$ such  that
             \begin{equation} \label{1.1}
\hbox{\rm Im}\,L_\la+ V=Y  \ \hbox{\rm for all }\  \la \in \Lambda.
\end{equation}
Because of \eqref{1.1}  the family  $E_{\la}=L_{\la}^{-1}(V),\, \la\in \Lambda, $ is 
a vector bundle $E.$ By definition, $\Ind  L = [E] \in \tilde{KO}(\Lambda ),$  where  $[E]$  denotes  the stable equivalence class of $E.$
 For families of Fredholm operators between complex Banach spaces the same construction produces an element  $\Ind L\in \tilde{K}(\Lambda).$

If $f(\la,x)=L_\la x,$ where  $\{L_\la \}_{\la \in \Lambda}$ is  a
family of linear Fredholm operators, then the set of singular points 
$ \Sigma(L)$  coincides  with the set of bifurcation points of $f$. By definition of   the index bundle,   $\Sigma(L)=B(f)$ is nonempty whenever $\Ind L\neq 0$  in $ \tilde{KO}(\Lambda).$ Hence, in the case of linear families,  bifurcation  is caused by the  nonvanishing of the index bundle.    However,  in order to detect the presence of bifurcation for a family of nonlinear Fredholm maps  $\Ind L$  is not sufficient,  and we have to resort to the image of $\Ind L$  by  the generalized $\J$-homomorphism  $\J\colon \tilde{KO}(\Lambda )\ra \J(\Lambda )$ \cite{[At-1]}.

Let us recall that two vector bundles $E,F$ are  fiberwise homotopy equivalent if there is a fiber preserving  homotopy equivalence  between the corresponding unit sphere bundles $S(E)$ and $S(F)$.  Moreover, $E$ and $F$ are said to be  {\it stably fiberwise homotopy equivalent }(shortly  sfh-equivalent) if they become  fiberwise homotopy equivalent after addition of  trivial bundles to both sides.  

Let $J(\LL)$ be  the   quotient group of $\tilde{KO}(\LL)$ by the subgroup generated by  elements of the form $[E] -[F]$ with $E$ sfh-equivalent to $F.$  The {\it generalized $J$-homomorphism}  $J\colon \tilde{KO}(\Lambda ) \ra J(\Lambda)$ is the projection to the quotient. By definition,  $J([E])$ vanishes in $J(\LL)$  if and only if $E$ is sfh-trivial.  The groups $J(\LL)$  were  introduced by Atiyah in \cite{[At-1]} who showed  that,  if $\LL$ is a finite $CW$-complex,   the group $J(\LL)$ is  finite.

  Let  $f\colon \Lambda \times O\ra Y$  be a continuous family of $C^1$-Fredholm maps (of index $0$)  parametrized by a finite connected $CW$-complex $\Lambda,$ such that  $f(\la,0)=0. $ The {\it index of bifurcation points} $\beta (f)\in \J(\LL)$ of the family  $f$ is defined by  $ \beta (f) =\J(\Ind L).$  
  
 Theorem $1.2.1$  in  \cite{[Pe]} states that, if 
 $ \beta (f) \neq 0$ in $J(\LL)$ and  $\Sigma(L)$ is a proper subset of $\Lambda,$ then  the family $f$ possesses at least one bifurcation point from the trivial branch.
  \vskip 3pt
  
  An $n$-dimensional vector bundle (n-plane bundle) is a vector bundle $\pi\colon  E \ra\LL$ such that $ \dim E_\la =n$ for all $\la\in \LL.$ Let  $Vect^n (\LL )$ be  the set of  all  isomorphism classes of $n$-plane bundles over $\LL.$ 
  
 If  $R$ is a ring, a characteristic class  $c \colon Vect^n (-) \ra H^{*}( -; R)$ is said to be of sfh-type (or spherical) if it depends only on the sfh-equivalence class of the vector bundle.
 
  Spherical characteristic classes  detect elements  with nontrivial  $J$-image.  Here we will consider only  Wu  classes  with values in the singular cohomology theory $H^{(2p-1)*}(-, \Z_p),$  if  $p$ is an odd prime, and  Stiefel-Whitney  classes, for  $p=2$.  
   
    Below we collect the needed properties of Wu classes. 
We will denote with   $X^*$  the Alexander one-point compactification of a locally compact space $X.$  $X^*$ is naturally a pointed space. 
  A proper map $f\colon X \ra Y$ extends uniquely to a map $\bar f \colon X^* \ra Y^*$ preserving base points.
  Moreover,  $(X\times Y)^*$  is homeomorphic to the wedge  
$$X^*\wedge Y^* = X^* \times Y^*/ ( X^*\times \{\infty\} \cup \{\infty\} \times Y^*).$$
This makes  the one point compactification into a  product preserving functor from the category of locally compact spaces  to the category of pointed compact spaces. 

 {\it Thom space} of an $n$-plane bundle  $\pi\colon  E \ra\LL $ is the one point compactification $E^*$ of its total space $E.$ Using the above homeomorphism  with $X=E$ and $Y $ a trivial $m$-plane bundle over a point  we conclude that the Thom space of   $ E\oplus \theta^m $  is the   $m$-th suspension of the Thom space of $E.$

Let $r_\la \colon E^*_\la \ra E^*$ be  the extension of the inclusion of the fiber  $E_\la$ into $E.$   An {\it orientation  (Thom) class}  for the vector bundle $E$ is an element $u \in \tilde H^n( E^*; \Z_p)$ such that, for all $\la,$  $r_{\la}^{*}(u)$ is a generator of $\tilde H^n( E_\la^*; \Z_p)  \simeq \Z_p.$

It is easy to see that every n-plane bundle admits an orientation over $\Z_2,$ and  that  a  bundle  is  orientable over $Z_p$ with $p>2$  if and only if it is orientable, i.e.,  it admits a reduction  of its structure group  to $SO(n).$   

The map  $d \colon E\ra \LL \times E$ defined by $d(v) =(\pi(v),v),$ being proper,  extends to a map $ 
\delta\colon E^* \ra \LL^* \wedge E^*.$  Composing the wedge product $$ \wedge \colon \tilde H^*(\LL^*; \Z_p)\times\tilde H^*( E^*; \Z_p) \ra \tilde H^*(\LL^*\wedge E^*; \Z_p) $$ with   $\delta^* \colon \tilde H^*(\LL^*\wedge E^*; \Z_p) \ra \tilde H^*( E^*; \Z_p)$ and using $H^*(\LL;\Z_p)\simeq \tilde H^*( \LL^{*}; \Z_p)$  we obtain a cup product 
$$ \cup \colon H^*(\LL; \Z_p)\times\tilde H^*( E^*; \Z_p)\ra \tilde H^*( E^*; \Z_p).$$  
  Thom's  isomorphism theorem   states that,  if $u\in \tilde H^n( E^*; \Z_p)$  is a Thom class for  $E,$  the homomorphism $$\Psi_u  \colon H^*(\LL; \Z_p) \ra \tilde H^{*+n}( E^*; \Z_p)$$ defined by  $\Psi_u(a)= a \cup u$ is an isomorphism.

 Let  $p$ be an odd prime. The {\it $k$-th Wu characteristic class}  $ q_k(E)\in H^{2(p-1)k}(\LL;\Z_p)$   of  an $n$-plane bundle $E$ orientable over $\Z_p$ is defined by  

\begin{equation}
\label{wu}
q_{k}(E)=\Psi_u^{-1} P^k u=\Psi_u^{-1} P^k \Psi_{u}(1),
\end{equation} 
where \[ P^k  \colon \tilde H^{n}(E^{*};\Z_p) \rightarrow \tilde H^{n+ 2(p-1)k}(E^{*};\Z_p)\]  is the $k$-th Steenrod reduced power \cite{[St-Ep]}. 

One of the consequences of the Thom isomorphism theorem is that $\tilde H^{*}(E^{*};\Z_p)$  is a free module over the ring $H^{*}(\LL;\Z_p) $  generated by $u$  via the cup product defined above.   

Any two Thom classes  $u,u' \in \tilde H^{n}(E^{*};\Z_p)$ are related by   $u= a\cup u'$ with $ a\in H^{0}(\LL;\Z_p) $ invertible. Since $P^{k}$ are module  homomorphisms substituting $u=a\cup u'$ in  \eqref{wu} we obtain that the  classes  $q_{k}(E)$ are  independent from the choice of the Thom  class $u.$ 
Since the suspension homomorphism commutes with $r^*_\la,$ 
from the characterizing property of  Thom's class  it follows   that the  $m$-th suspension  $ u'=\sigma^{m}(u)$  of a Thom class $u$ of $E$  is a Thom class for $E \oplus\theta^{m}.$ Moreover, since  the cup product verifies $ a \cup \sigma u = \sigma( a \cup u),$ we have   $\Psi_{u'}=\sigma^{m}\Psi_{u}.$  On the other hand also  $P^{k} $ commute with the suspension.  Hence,  we get
 \begin{equation}\label{stab}
q_{k}(E\oplus \theta^{m}) =\Psi^{-1}_{u'} P^{k} u' = \Psi_{u}^{-1} \sigma^{-m}P^k \sigma^{m}u =\Psi^{-1}_{u} P^{k}u= q_{k}(E). 
\end{equation}

Thus  $q_{k}$ depends only on the stable equivalence class of $E,$ and hence we have a well defined natural transformation $ q_{k }\colon \tilde{KSO} (-)\ra  H^{ 2(p-1)k}(-;\Z_p),$  where $\tilde{KSO}(-)$  is the ring of stable equivalence classes of orientable bundles.  As a matter of fact, the classes  $q_k$ can be  defined for all elements of $\tilde{KO}(-),$ but we will not use this here.

On the other hand, a  fiberwise homotopy equivalence  $h \colon S(E) \rightarrow S(F)$ by radial extension produces a proper homotopy equivalence between the total spaces of $E$ and $F$  and hence a base point preserving homotopy equivalence $\bar h\colon E^{*} \rightarrow F^{*}.$ This later restricts to a homotopy equivalence $\bar h_\la \colon E_\la^{*} \rightarrow F_\la^{*}.$ It follows from this  that,  if $v$ is an orientation for $F,$ then $u={\bar h}^{*}(v)$ is an orientation for $E$ and  moreover $\Psi_{u}={\bar h}^{*}\Psi_{v}.$  Substituting  in \eqref{wu}  we get $q_k(E)=q_k(F).$ Thus   $q_{k}$ depends only on the fiberwise  preserving homotopy class of the sphere bundle $S(E)$ and hence  $q_{k}\colon \tilde{KSO}(-)\ra H^{ 2(p-1)k}(-;\Z_p) $  factorizes through the functor $J(-).$
 
   The same holds for the Stiefel-Whitney classes $\omega_k \in H^{k}(-;\Z_2)$  since they are  constructed from the Thom class of $E$  using  Steenrod squares $Sq^{k}$  in \eqref{wu}. 
    
With this at hand we can state the main result of this section. 

 Let us first  recall that the covering dimension of a topological space $X$ is defined to be the minimum value of $n$  such that every open cover of $X$ has an open refinement  for which  no point is included in more than $n+1$ elements. By a well known characterization due to Hurewicz,  the  topological dimension of a compact space $X$ is at least $n$ if, for some closed subset $C$ of $X,$   the Alexander-Spanier cohomology  $H^{n}(X,C)\neq 0.$

  \begin{theorem}\label{th:1.2} 
  
 Let  $ \Lambda$ be a compact  connected topological manifold and let \\ $f\colon\Lambda \times O\rightarrow Y$  be a continuous family of $C^1$-Fredholm maps  verifying  $f(\la,0)=0$ and such that $\Sigma(L)$ is a proper subset of $\LL.$ 
 
\begin{itemize}
\item[i)] \  If $\Lambda$ and $\Ind L$ are orientable and,  for some odd prime  $p,$ there is a  $k\geq 1$  such that   $q_k(\Ind L)\neq 0 $ in $H^{2(p-1)k}(\Lambda;\Z_p),$ then the Lebesgue covering dimension of the set $ B(f)$ is at least $\dim \LL - 2(p-1)k.$

\item[ii)] \  If  $\omega_k(\Ind L) \neq 0$  in $H^k(\Lambda;\Z_2) $ for some  $k\geq 1,$ then the
dimension of  $B(f)$ is at least $\dim\LL -k.$
\end{itemize}
Moreover,  either the set $B(f)$ disconnects $ \Lambda$ or  it cannot  be  deformed in $ \Lambda$ into a point.
\end{theorem}

\proof  (see also \cite{[Fi-Ma-Pe]}, \cite{[Ba]})  We will denote  with $\bar H^{*}(X;\Z_p) $ the  Alexander-Spanier cohomology of $X$ with $\Z_p$  coefficients.  It is well known that  $\bar H^{*}(X;\Z_p)$ coincides with the singular  cohomology of $X$ when $X$ is a manifold.

 Let $B=B(f).$  If $\dim \LL =m$  and $\Lambda-B$ is not connected then the covering dimension of $B$  must be
 at least $m-1,$ since sets of smaller dimension cannot disconnect  $\Lambda.$ Hence, in this case, the theorem is proved. 
 
 From now we assume  that  $\Lambda -B$ is connected.  Let $\alpha \in H_{2(p-1)k}(\Lambda ;{\Z_p})$ be any homology class such that the Kroenecker pairing  $<q_k(\Ind L);\alpha>\neq 0$ and let $\eta\in  H^{m-2(p-1)k}(\Lambda; \Z_p) $ be the  Poincar\'e dual of $\alpha.$ 
 
  Let $i\colon B\ra \LL$ be the inclusion and let $\zeta = i^*(\eta) \in \bar H^{*}(B;\Z_p) $ be the restriction of $\eta $ to 
  $B.$  If we  can show that  $ \zeta \neq 0 $ in $ \bar H^{*}(B;\Z_p),$  then the theorem is proved. Indeed, $\zeta = i^*(\eta)$ is  an obstruction to the deformation of the  subspace $B$ to a point and,  by cohomological characterization of the covering dimension, $\dim B$ must be at least  $m-2(p-1)k$. 
  
In order to show that  $\zeta\neq 0$  let us consider the following 
commutative diagram

\begin{equation}
\label{eq:1.4}
\begin{matrix} \ &\ &\ & i^* & \ \cr
\ &\ & H^{m-2(p-1)k}(\Lambda ;{\mathbb Z_p}) & \longrightarrow & \bar H^{m-2(p-1)k}(B;{\mathbb Z_p}) \cr
\ &\ & \Big\uparrow & \  &\Big\uparrow \cr
 H_{2(p-1)k}(\Lambda  - B;{\mathbb Z_p})& \longrightarrow &  H_{2(p-1)k}(\Lambda ;{\mathbb Z_p})
 & \longrightarrow & H_{2(p-1)k}( \Lambda , \Lambda -B;{\mathbb Z_p})\cr
\ & j_* &\ &\pi_* &\ \cr 
\end{matrix}
\end{equation}
 where  the vertical arrows are the  Poincar\'e duality isomorphisms and the bottom sequence is the homology sequence of a pair.

By commutativity, $\zeta$ is dual to $\pi_*(\alpha).$ Hence, it is enough to show that the homology class $\pi_*(\alpha)$ does not vanish. 

If so, by exactness, 
$\alpha= j_*(\beta)$ for some $\beta\in  H_{2(p-1)k}(\Lambda -B;{\mathbb Z_p}).$ Since singular homology has compact supports, there  exists a  finite connected  $CW$-complex  $P$ and a map   $h \colon P\ra \Lambda-B $ such that $\beta = h_*(\delta )$  for some $\delta \in H_{2(p-1)k}(P; {\mathbb Z_p})$ (for this  it is enough to take as $P$ any closed  connected  polyhedral neighborhood of the support of a singular cochain representing $\beta).$

Since  $\Lambda-B$ is connected,  we can  assume without loss of generality  that  some point  $ \la_0 \in  \LL -\Sigma(L)$ belongs to the image of $ h.$    Let us consider now $\bar h=j h$ and the family $\bar f\colon P\times X\ra Y$ defined by $\bar f (p,x)= f(\bar h (p),x).$ 

The linearization at the trivial branch of $\bar f$ is $\bar L =  L  \bar h.$  Since   $\la_0=h(p_0)$ is not  a singular of $L,$  the set $\Sigma (\bar L)$ is a proper subset of $P.$  On the other hand, by construction, $\bar h$ sends bifurcation points of $\bar f$ into bifurcation points of $f,$ and since $ \bar h(P) \cap  B =\emptyset $, the family $\bar f$ has no bifurcation points.

By   \cite[ Theorem $1.2.1$ ]{[Pe]},    $J(\Ind \bar L) =0 $  and hence all characteristic classes $q_{k}$  of  $\Ind \bar L$ must vanish. But  $ \Ind \bar L = \bar{h}^*(\Ind L),$  and  by naturality of characteristic classes 
$q_k (\Ind (\bar L ))= h^* j^*q_k( \Ind L).$

Hence,
$$0 =<h^* j^*\bigl[q_k( \Ind L)\bigr]; \delta>=  <q_k(\Ind L);j_*(\beta)> = <q_k( \Ind L);\alpha>,$$
which contradicts the choice of $\alpha.$

The  proof of $ii)$ is similar. 

\vskip 5pt

\section{Proof of the main theorem}

Denoting  by  $H^s$ the Sobolev-Hardy spaces,  the map  $$(\cf , \cg) \colon \Lambda  \times  C^\infty(\Omega;\R^m ) \rightarrow  C^\infty(\Omega;\R^m ) \times  
\prod _{i=1} ^r C^\infty(\partial \Omega;\R )$$ defined by \eqref{cinf} 
extends to a smooth map 
\begin{equation}
\label{nbvp1}
 h=(f,g) \colon \Lambda \times  H^{k+s}(\Omega;\R^m ) \rightarrow  H^{s}(\Omega;\R^m ) \times  H^{+}(\partial \Omega;\R^r), 
\end{equation}
where,  by definition, 
$H^{+}(\partial \Omega;\R^r) = \prod _{i=1} ^r H^{k+s-k_i-1/2}(\partial \Omega;\R ).$

   By our assumptions,   $u\equiv 0$ is a solution of $h_\la(u)=0$  for all $\la\in \Lambda.$  Hence $\LL\times\{0\}$ is   a trivial branch for $h.$   The Frechet derivative  $Dh_\la$ at $u=0$  is the operator $(L_\la, B_\la)$ induced on Hardy-Sobolev spaces  by $(\cL_{\la}, \cB_\la).$  
     
    It is shown in  \cite[Proposition 5.2.1]{[Pe]}  that, under the  hypothesis of Theorem \ref{th:0},  we can find a  neighborhood  $O$  of $0$ in $\hk$  such that 
 \[h \colon \Lambda \times O \rightarrow  H^{s}(\Omega;\R^m ) \times H^{+}(\partial \Omega;\R^r)\]   is a smooth parametrized family of  Fredholm maps of index $0.$  By \cite[Proposition 5.2.2]{[Pe]}, the set of bifurcation points of   \eqref{bvp1} coincides with the set $Bif(h)$ of bifurcation points of the map $h$.
Moreover, denoting with $(L,B) $ the linearization of $h$  along  the trivial branch, we have that $\nu \notin \Sigma(L,B)$. 

 If $\Ind (L,B)$ is nonorientable,   $\omega_1\Ind (L,B)\neq 0.$  By  assertion $ii)$ of  Theorem \ref{th:1.2},  $\dim B \geq  d-1$  and $Bif(h)$ carries a nontrivial class of positive degree in cohomology with $Z_2$ coefficients.   Hence,   in this case, the conclusions of Theorem \ref{th:0} hold regardless any condition on $d(\si,\Sigma).$ 

If  $\Ind (L,B)$ is orientable,  the proof of Theorem \ref{th:0}  is obtained by relating the degree  $\deg( \sigma;\Sigma) $ defined in \eqref{fed} with the evaluation of the first Wu class of $\Ind(L,B)$ on the spherical homology class  $[\Sigma].$    Only the first Wu class is of interest for us  because, as we will see, all characteristic numbers obtained in this way from higher Wu classes vanish. 

The  relation between  characteristic classes of the index bundle  and the degree $\deg(\sigma;\Sigma)$ comes from the family version of the Agranovich reduction and the Atiyah-Singer theorem.

  Let $\sigma \colon \LL\times (R^{2n} - K\times\{0\}) \to GL(m;\C)$ be  the map defined by \eqref{symbo}.  In   \cite{[Pe]} we have constructed a smooth family  $ \cS \colon \LL \ra \pc $  of elliptic  pseudo-differential operators of order zero on  $\R^{n}$  such that the principal symbol of $\cS_{\la}$ coincides with $\sigma_{\la}.$ The Agranovich reduction relates  $\cS$ with the family $(\cL,\cB)$ considered as a family of differential operators with complex coefficients. 
  
 More precisely,  denoting  by  $S_{\la}\colon  H^{s}(\Omega;\C^m )\ra  H^{s}(\Omega;\C^m )$  the  operator induced by $\cS_{\la}$  on Hardy-Sobolev spaces and  with  $(L^{c},B^{c})$  the complexification of $(L,B),$ Theorem $4.1.1$ of \cite{[Pe]} states that  in $\tilde K(\Lambda)$
   \begin{equation}
\label{agranred}
\Ind ( L^{c},B^{c})= \Ind S.
\end{equation}

Since  the index bundle of the family $(L^{c},B^{c})$  is the complexification of the (real) index bundle $\Ind (L,B),$ we have:  
  \begin{equation}\label{ag}
c(\Ind (L,B))= \Ind S,
\end{equation}
where  $ c \colon \tilde{KO} (\LL) \rightarrow \tilde K(\LL)$ is the complexification homomorphism.

 Using the Chern-Weil theory of characteristic classes of smooth vector bundles over smooth manifolds,  in \cite{[Fe-1]}  Fedosov obtained  an explicit expression for the differential  form representing  Chern character  $\ch(\Ind S)$  in de Rham cohomology with complex coefficients.  He   showed that 
$\ch (\Ind S) $ is the cohomology class of the form
\begin{equation}\label{fedfor}
-\sum_{j=n}^\infty\displaystyle {\frac{(j-1)!} {(2\pi i)^j(2j-1)!}} \oint _{S^{2n-1}} tr(\sigma^{-1}d\sigma)^{2j-1},
\end{equation} 
where  $\oint$ denotes  the {\it integration along the fiber} (see \cite{[Pe]} Appendix C) and  $S^{2n-1}$ is the boundary of  a ball in $\R^{2n}$ containing the set  $\{(x,\eta) / \det \sigma (\la,x,\eta)=0\}.$ 

Using Fedosov's  formula  we will show that, under the hypothesis of theorem \ref{th:0},  $q_{1} (\Ind(L,B))\neq 0$ in $H^{2(p-1)}(\LL;\Z_p).$  Then theorem \ref{th:0}  will follow immediately  from   Theorem \ref{th:1.2}  i)   and  \cite[Proposition 5.2.1]{[Pe]}, which shows that the set $B(h)$ of bifurcation points of the map  $h$ coincides with the set  $B$ of bifurcation points for  classical solutions of the system \eqref{bvp1}.  

  The  rest of this  section is devoted to show  that $q_{1} (\Ind(L,B))\neq 0.$  For this, we will consider the restriction of the family  $h$ to $\Sigma \times  H^{k+s}(\Omega;\R^m ).$  
  
   More precisely, if $q=2(p-1)$ and  $e\colon S^{q} \ra \Lambda$ is  an orientation preserving  embedding  with $\im e=\Sigma,$ let $$\bar h \colon S^{q} \times  H^{k+s}(\Omega;\R^m ) \rightarrow  H^{s}(\Omega;\R^m ) \times H^{+}(\partial \Omega;\R^r) $$  be defined  by $\bar h(\alpha,u)=h(e(\alpha),u).$ 
  
   The family $\bar h$ is the nonlinear Fredholm map induced in functional spaces by the pullback of the problem \eqref{bvp1}  to the sphere $S^{q}.$   The linearization  of $\bar h $ 
  along the trivial branch  is  the family  $(\bar L,\bar B)= (L,B)\circ e.$

Let  $ \bar\sigma (\alpha,x ,\eta) =\sigma(e(\alpha) ,x,\eta) $ and  $\bar S =S\circ  e.$ Then  $\bar S$ is induced in functional spaces by the family of pseudo-differential operators $\bar \cS$ with principal symbol $\bar \sigma.$

By the previous discussion, from \eqref{ag} we get 
\begin{equation}\label{pseudo}
 c(\Ind (\bar L,\bar B))= \Ind \bar S
 \end{equation}
 Let us show that, for $q=2(p-1),$ the Kroenecker pairing \[ <q_{1}(\Ind(L,B)); e_{*}([S^{q}] > \neq 0.\] 

By naturality of characteristic classes  and the index bundle we have 
\begin{equation}\label{proof}
<q_{1}(\Ind(L,B)); e_{*}([S^{q}]) >= <q_{1}(\Ind(\bar L,\bar B)); [S^{q}]>
\end{equation}
In order to compute the right hand side let us recall that, putting $r=\frac12 (p-1),$  the Wu class  $q_{k}(E)$  of an n-plane bundle over $X$  can be written as a polynomial  $K_{rk}(p_{1},\cdots,p_{rk})$ in Pontriagin classes  reduced mod $p.$  The polynomials $K_{rk}$  are  with $\Z_{p}$ coefficients  and belong to a multiplicative sequence associated to the function $\phi(t)=1+t^{r}$ \cite{[Mi]}.

To shorten notations,  let  $\eta = \Ind(\bar L,\bar B)$ and  $\eta^c = c(\Ind(\bar L,\bar B)).$ 

For vector bundles over $ S^{4r}$ we have  $p_{i}(\eta) =0$ for $i<r$  and  $q_{1}(\eta)= \pm r p_{r}(\eta)$ reduced mod $p$.  Indeed, it follows from Lemma 1.4.1  in \cite{[Hi]}  and  Newton's identity relating power sums to elementary symmetric functions, that,  over the integers, the coefficient of the integral Pontriagin class $p_r$ in $K_r$ is given by $s_r(0,\dots,0,1) =\pm r.$
 
  By Bott's  integrality  theorem \cite[section 18.9]{[Hu]},  the Chern character  $ch(\eta^c)=ch_{2r}(\eta^c)$ is an integral class. Moreover, using   \eqref{pseudo}, for the integral  Pontriagin class  $p_{r}$ it holds that
\begin{equation}
p_{r}(\eta)=(-1)^{r} 
c_{2r}(\eta^c)= \pm(2r-1)!ch_{2r}(\Ind \bar S),
\end{equation}
which gives \begin{equation}\label{wudeg}
<q_{1}(\Ind(\bar L,\bar B)); [S^{4r}]>= \pm r(2r-1)! <ch_{2r}( \Ind \bar S); [S^{4r}]>  \mod p.
\end{equation}
Here we denote in the same way the fundamental class in homology with coefficients in $\Z$ and $\Z_{p}.$

By  Fedosov's formula \eqref{fedfor},  the differential form representing
$ch_{2r}(\Ind \bar S)$
in de Rham cohomology is
\[  \Omega ={\frac{(n+2r-1)!} {(2\pi i)^(n+2r)(2n+4r-1)!}} \oint _{S^{2n-1}} tr(\bar \sigma^{-1}d\bar \sigma)^{2n+4r-1}.\]

Since  the cohomology class of $\Omega$ belongs to $H^{{4r}}(S^{4r};\Z),$ we have that 
 
 \begin{eqnarray*}
  <ch_{2r}(S); [S^{4r}]> = {\frac{(n+2r-1)!} {(2\pi i)^(n+2r)(2n+4r-1)!}}  \int_{ S^{4r}} \oint _{S^{2n-1}} tr(\bar \sigma^{-1}d\bar \sigma)^{2n+4r-1} =\\ =  {\frac{(n+2r-1)!} {(2\pi i)^(n+2r)(2n+4r-1)!}} \int_{S^{4r}\times S^{2n-1} } tr(\bar \sigma^{-1}d\bar \sigma)^{2n+4r-1}
\end{eqnarray*} is an integer.  

The last term of the above expression coincides with  $\deg(\sigma;\Sigma)$ defined in \eqref{fed}. Thus, from \eqref{wudeg} we obtain
 \[ <q_{1}(\Ind(\bar L,\bar B)); [S^{4r}]>= \pm r(2r-1)! \deg(\sigma;\Sigma)  \mod p.\]  
 Since $ r(2r-1)!$  is not divisible by $p,$  by  \eqref{proof}   $q_{1}(\Ind(L,B))\neq 0$ whenever $\deg(\sigma;\Sigma)$ is not divisible by $p.$  This concludes the proof of the theorem.  \qed

\br {\rm  Notice that  the above calculation gives 
$$<q_{k}(\Ind(\bar L,\bar B)); [S^{4r}]>=0,  \hfor  k>1.$$} \er 

\vskip 5pt 

\section{Other results}
In this section,  we will obtain sufficient conditions for bifurcation of solutions for  a particular class of nonlinear elliptic boundary value problems   \eqref{bvp1}  on dimensions not covered by Theorem \ref{th:0}.  

 More precisely, 
we substitute $m$  with  $ m'=2m$ and $r$ with $r'=2r$  in \eqref{bvp1} and we denote with $(\cL',\cB')$ the linearization  of \eqref{bvp1} along the trivial branch in both hypotheses  
$H_2$ and $H_3.$ But  instead of  $H_1$ we assume:

$H'_1)$ - The principal part of the linearization $(\cL',\cB')$ at the trivial branch is obtained from a family of elliptic boundary value problems for linear partial differential operators with complex coefficients 
 \[(\cL,\cB)\colon\R^q\times C^\infty( \Omega;\C^m) \ra C^\infty( \Omega;\C^m) \times\ C^\infty(\partial\Omega;\C^r)\] 
by forgetting the complex structure. 

Notice that   $H'_1$ is verified by semilinear equations  arising as real  lower-order perturbations of  families of linear elliptic boundary value problems with complex coefficients.  It  also holds in the case of  quasilinear elliptic problems for functions with values in $\C^m$  of the form    \be\label{qlfm} \cf(\la,x,u, D^k u) = \sum_{|\alpha|=k}  a_\alpha (\la, x, u, D^{k-1}u)D^\alpha u +  l.o.t.,\ee   where  $ a_\alpha \in \C^{m\times m},$  and  the boundary map $g$ defined in the same way.   In this case the assumption $H'_1$ is verified  because, for maps $\cf_\la  $  defined by \eqref{qlfm},  the principal part of the linearization at $0$  is given by  $\cL_\la  v = \sum_{|\alpha|=k}  a_\alpha (\la, x, 0)D^\alpha v,$  and similarly for $\cg.$

Assuming $H'_1,$ we will  use the principal symbol  $p(\la,\xi)$ of the complex differential operator $\cL$ to  define  \begin{equation} \label{symbo'}  \sigma \colon \LL\times (\R^{2n} - K\times\{0\}) \rightarrow GL(m;\C)\end{equation} in the same way as  in \eqref{symbo}. Moreover, given  a $q-$sphere $\Sigma$ embedded in $\Lambda,$  we define   $\deg(\sigma;\Sigma)$ by the equation \eqref{fed}.
 
 \begin{theorem}\label{th:0'}
Let the boundary value problem \eqref{bvp1} verify $H'_1,H_2$ and $H_3.$ If, for some sphere $\Sigma $ of dimension $q=2,4$ embedded in $\Lambda,$ the number  $\deg(\sigma;\Sigma )$ defined by \eqref{fed} is odd,  then
 the Lebesgue covering dimension of the set $B$ of 
all bifurcation points of \eqref{bvp1} is at least $d-q$ 
and the set $B(f)$ either disconnects $ \Lambda$ or is not contractible in  $ \Lambda$ to a point. 
\end{theorem}

\proof: A  complex Fredholm operator of index $0$ is still Fredholm, with the same index, when viewed as a real operator. Hence, from $H'_1$ it follows that  the family  $(L'B')$  induced by $(\cL',\cB')$ is a family of Fredholm operators of index $0$ and moreover  $(L'_\nu, B'_\nu)$ is invertible. 

In the same way as in the previous section, we can find a neighborhood  $O$ of $u\equiv 0$ such that      $ h= (f,g) \colon \Lambda \times O \rightarrow  H^{s}(\Omega;\R^{m'} ) \times H^{+}(\partial \Omega;\R^{r'})$  is a smooth parametrized family of  Fredholm maps of index $0$ whose linearization at the trivial branch  is the family $(L',B').$ 

If $\cS$ is the family of pseudo-differential operators associated to $\sigma,$ then,  by Theorem $4.1.1$ of \cite{[Pe]},
$\Ind(L,B)=\Ind S.$

Let  $r \colon  \tilde K(\Lambda)  \to  \tilde{KO}(\Lambda)$ be the  homomorphism obtained  by forgetting the complex structure.  We will  compute the Stiefel-Whitney classes  of   $ r \left(\Ind (S)\right)$ by restricting them  to $\Sigma.$
 As before, let $e$ be an embedding of the sphere $S^q, q=2s,$ whose  image is  $\Sigma.$
For any element   $\eta\in \tilde K(S^{2s}),$  the Chern character $ch (\eta)=ch_s(\eta)$ is an integral class. Moreover  the Chern class   $c_s (\eta) = \pm (s-1)! ch_s(\eta).$

 The Stiefel-Whitney class  $\omega_q(\eta)$ is the  image of $c_s$ under  the change of coefficients 
 $ \rho \colon H^q(S^q;\Z)  \ra H^q (S^q;\Z_2)$ and therefore,
   $<\omega_q r(\eta); [S^q]>$ is the mod $2$ reduction of $(s-1)! <ch_s(\eta);[S^q]>.$ 
  Thus, if  $s=1,2$ and $ <ch_s(\eta);[S^q]>$ is odd, $\omega_q r(\eta)\neq 0$  in $ H^q(S^q;\Z_2).$    
 
 Using Fedosov's  formula for the Chern character of the index bundle as before, we obtain $<ch_s(e^*(\Ind S) );[S^q]>=\deg(\sigma;\Sigma).$
  Taking   $\eta =e^*( \Ind S)$ in the above discussion,  if  $q=2,4$ and  $d(\sigma;\Sigma)$ is odd,  then   $e^*\omega_q r( \Ind S) \neq 0 \ \text{in}\  H^q (S^q;\Z_2)$ and hence
  $\omega_q r \left(\Ind (L,B)\right)= \omega_q  r (\Ind S)\neq 0$ in $H^q (\Lambda ;\Z_2)$ as well. 
Thus  $B(h)$ verifies the conclusion $ii)$ of  theorem  \ref{th:1.2} and the bootstrap  \cite[Proposition 5.2.1]{[Pe]}  concludes the proof.
   \qed

\end{document}